\documentclass[12pt]{article}
\usepackage[utf8]{inputenc}

\usepackage[top=1in, bottom=1in, left=1in, right=1in]{geometry}
\usepackage{graphicx} 
\usepackage{subfig}
\usepackage{amsthm}
\usepackage[authoryear]{natbib}
\usepackage{ulem}
\usepackage{color}
\usepackage{lineno}  
\usepackage{hyperref}
\usepackage{amsmath}
\usepackage{amssymb}
\usepackage{indentfirst} 

\numberwithin{equation}{section}

\newcommand\keywords[1]{\textbf{Keywords}: #1}

\title{Linear stability of inner case of double averaged \\ spatial restricted elliptic three-body problem}
\author{\small{Xiumin Huang$^{1,2}$, Yan Luo$^3$, Kaicheng Sheng$^4$\footnote{Corresponding author. E-mail addresses:  {\it k.sheng@sdu.edu.cn} (K.Sheng)}}\,, Yiru Ye$^5$ \\
\small{$^1$ School of Astronomy and Space Science, Nanjing University, Nanjing, China }\\
\small{$^2$ Key Laboratory of Modern Astronomy and Astrophysics in Ministry of Education, } \\ \small{Nanjing University, Nanjing, China} \\
\small{$^3$ Research Centre for Mathematics and Interdisciplinary Sciences,} \\ \small{ Shandong University, Qingdao, China }\\
\small{$^4$ School of Mathematics, Shandong University, Jinan, China }\\ 
\small{$^5$ School of Mathematics and Physics, Xi'an Jiaotong-Liverpool University, Suzhou, China}
}

\date{\today}

\begin{document}


\maketitle

\begin{abstract}
We study the secular effects in the motion of an asteroid with negligible mass in a spatial restricted elliptic three-body problem with arbitrary inclination. Averaging over mean anomalies of the asteroid and the planet are applied to obtain the double averaged Hamiltonian system. It admits a two-parameter family of orbits (solutions) corresponding to the motion of the third body in the plane of primaries’ motion. The aim of our investigation is to analyze the stability of these orbits in inner case. We show that these orbits are stable in the linear approximation with spatial perturbation and give descriptions of linear stability with respect to the eccentricity and argument of periapsis of asteroid. Numerical simulations of different types of orbits are performed as well. 

\bigskip

\noindent\keywords{Three-body problem, Averaging procedure, Lidov-Kozai effect, Linear Stability.}
\end{abstract}



\section{Introduction}

Averaging method has significant applications in celestial mechanics. This paper focuses on the double averaged spatial elliptic restricted three-body problem (ER3BP) involving a star, a planet, and an asteroid. The study is under the assumptions that the mass of the planet is much smaller than the mass of the star and the distance between the asteroid and the star is much smaller than the distance between the planet and the star (referred to as inner case or the case of a distant perturber). The averaging over fast phases corresponds to the motions of the star-planet system and the asteroid.

If the primary bodies (i.e., the star and the planet) move around their barycenter in circular orbits, it is called circular problem. In this case, the double averaged problem is integrable ~\citep{moiseev}. Under the inner case, the problem has been analyzed by~\citet{zeipel},~\citet{lidov} and~\citet{kozai}. Complete numerical simulation of bifurcations in the double averaged restricted circular three-body problem is accomplished by~\citet{vashkovyak}. 

When the primary bodies move in elliptical orbits, the double averaged three-body problem becomes non-integrable and shows chaotic dynamics in its phase space~\citep{Naoz2016}. However, for inner case or the distance between the asteroid and the star is significantly larger than the distance between the planet and the star (in the case of double star), integrable approximations of the elliptic double averaged problem can be performed~\citep{lidov,zig}. These approximations involve truncating the expansion of the disturbing function with respect to the ratio of the mentioned distances at the principal term. In inner case, this approximation does not differ qualitatively from the case that the primaries move in circular orbits.~\citep{lidov}. Consideration of more accurate models can be explored in works such as~\citet{Katz2011},~\citet{Lithwick2011} and~\citet{Sidorenko2018}. For the complete study of the double star case, refer to the work by~\citet{zig}.

Trajectories corresponding to the motion of the asteroid in planar orbit (i.e., in orbit which is in the plane of primaries' motion) is an invariant manifold of the restricted three-body problem. In circular problem, the planar orbits undergo a secular evolution characterized by uniform precession. However, when the primaries move in elliptical orbits (elliptic problem), the evolution becomes more intricate. For inner case of elliptic problem, the evolution of the planar orbits has been investigated analytically by~\citet{akse}. A completed numerical analysis, including all possible situations for the evolution of planar orbits in the ER3BP, can be found in the work by~\citet{vash_planar}.

The discussion of the planar orbits naturally implies the important question of their stability with respect to spatial perturbations.
Such a stability of planar orbits in double averaged circular problem was established by \citet{neish}.  His study indicates that planar orbits also remain stable in the linear approximation of the double averaged elliptic problem with a sufficiently small eccentricity of the perturber's orbit. For the inner case and the case of double star in spatial perturbed problem, the planar orbits are stable within the scope of the double averaged problem, regardless of the eccentricity of the perturber's orbit~\citep{lidov,zig}.  The stability of planar orbits which are equilibria of double averaged elliptic problem was studied by \citet{NSS} in spatial perturbed problem. A small inclination was assumed to be the perturbation in order to simplify the investigation of the stability. The considered equilibria are stable in the linear approximation. It turns out that the resonance 2:1 between frequencies of oscillations in the plane of the primaries' motion and across this plane leads to instability (at least for the limiting case of large ratio in semi-major axes of the asteroid and the planet). 

The stability of planar orbits with respect to arbitrary inclinations instead of spatial perturbations could be rather complicated. Orbital Flip cases are studied by~\citet{Katz2011},~\citet{Sidorenko2018} and~\citet{Lei}.

It is natural to study the problem in inner case first. The goal of this paper is to investigate the stability of planar orbits in inner case of the double averaged ER3BP with respect to a large variation of the inclination. {Similar to the work in Neishtadt-Sheng-Sidorenko (2021), we consider the equilibria as well as the periodic orbits around the equilibria in planar problem and then give spatial perturbation to them. We aim to determine if the equilibria and periodic orbits are still stable under the spatial perturbation. The cases of small spatial perturbation (inclination $i$ can be expanded) and large spatial perturbation (inclination $i$ can not be expanded) are discussed. Hamiltonian system of planar problem has 2-degree of freedom (2-DOF), with spatial perturbation, it changes to 3-degree of freedom (3-DOF), linear stability of the equilibria and periodic orbits of the asteroid in double averaged planar problem is studied.} Description of linear stability with respect to the eccentricity and argument of periapsis of the asteroid is given in a general case. 

\section{Statement of the problem}

\begin{figure}[htbp] 
\centering
\includegraphics[width=0.5\textwidth]{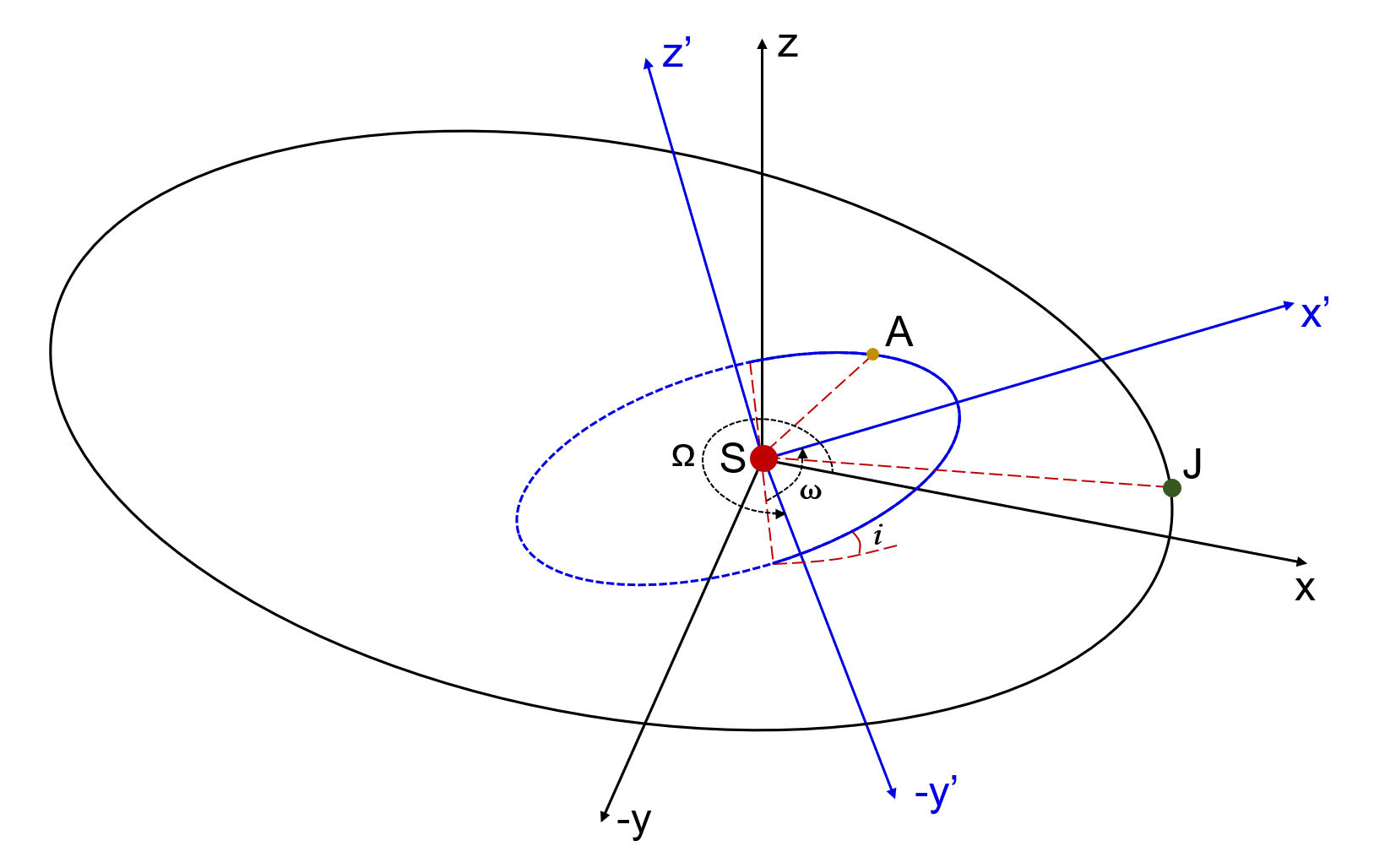} 
\caption{Coordinate frames.}
\label{fig:tu77} 
\end{figure} 


We consider a spatial elliptic restricted three-body problem involving a star $S$, a planet $J$ and an asteroid $A$ \citep{BC}. Similar to the frame in \citet{NSS}, we set the position of the star as the origin $O$ of the Cartesian coordinate system $Oxyz$, and the plane $Oxy$ of the system is determined by the motion of the star and the planet. Thus the coordinates of the planet and the asteroid are $(x_{J},y_{J},0)$ and $(x,y,z)$, respectively. The Cartesian coordinate frame $Ox'y'z'$ is a rotating frame of $Oxyz$, where the plane $Ox'y'$ is the osculating plane of the orbit of the asteroid,  then $(x', y',0)$ are coordinates of the asteroid in the rotating frame. We use the standard osculating elements $a$, $l$, $e$, $\omega$, $i$, $\Omega$ to describe the orbit of the asteroid, which represent the semi-major axis, mean anomaly, eccentricity, argument of periapsis, inclination, and longitude of the ascending node, respectively. See Fig. \ref{fig:tu77}.  {We study the spatial problem where the inclination $i$ can be considered as a perturbation of the planar problem.} We have

\begin{equation}
\label{coord_transform}
\begin{array}{ccl}
x &=&\left( \cos\Omega
\cos\omega -\cos  i\sin\Omega
  \sin\omega  \right) x'+\left( -\cos\Omega\sin  
\omega  -\cos i\sin\Omega \cos
\omega  \right) y',
\\ 
 y &=&\left( \sin 
\Omega  \cos\omega +\cos i\cos\Omega
 \sin\omega  \right) x' +\left( -\sin\Omega\sin
\omega +\cos  i\cos\Omega
  \cos\omega  \right) y' ,
\\ 
 z &=&\left(\sin  i
  \sin\omega\right) x'+ \left(\sin i \cos\omega\right) y'. 
\end{array} 
\end{equation}

It follows from \citet{shev}, the planet moves in a prescribed elliptic orbit:
\begin{eqnarray*}
\label{coord_J}
 {{x}_{J}}={{a}_{J}}\left( \cos{{E}_{J}} -{{e}_{J}} \right),\
  {{y}_{J}}={{a}_{J}}\sqrt{1-{e}_{J}^{2}}\, \sin  {{E}_{J}} ,\quad
  l_J=E_J-e_J\sin E_J. 
\end{eqnarray*}
Here $a_J$, $e_J$, $E_J$, and $l_J$ are the semi-major axis, the eccentricity,  the eccentric anomaly, and the mean anomaly of the planet's orbit.  We put $a_J=1$ for convenience in the following.

\section{Hamiltonian of the system}
We introduce the canonical Delaunay elements $l$, $g$, $h$, $L$, $G$, $H$, where $l$, $g=\omega$, $h=\Omega$ are the mean anomaly, argument of pericenter, ascending node of the asteroid, respectively. $L=\sqrt{(1-\mu)a}$ corresponds to Keplerian energy, $G=L\sqrt{1-e^{2}}$ is the total angular momentum and $H=G\cos i$ is component of angular momentum perpendicular to the equator \citep{BC}. 

We consider canonical Poincar\'e variables $p_1$, $q_1$, $p_2$, $q_2$, $p_3$, $q_3$ to demonstrate the dynamics of the asteroid:

\begin{equation}
\label{P_elements}
\begin{array}{ll}
  {{p}_{1}}=L, & {{q}_{1}}=l+g+h,\\ 
  {{p}_{2}}=\sqrt{2\left( L-G \right)}\cos\left( g+h \right),\quad &
  {{q}_{2}}=-\sqrt{2\left( L-G \right)}\sin \left( g+h \right),
  \\ 
 {{p}_{3}}=\sqrt{2\left( G-H \right)}\cos  h ,& 
  {{q}_{3}}=-\sqrt{2\left( G-H \right)}\sin  h.
\end{array}  
\end{equation}

Let the mass of the planet be $\mu$ and the mass of the star be $1-\mu$. Thus the sum of the mass is $1$. The Hamiltonian of the asteroid is
\begin{equation}
\label{F}
F=-{\frac {\left( 1-\mu \right) ^{2}}{2{L}^{2}}}+\mu U, 
\end{equation}
where
\begin{equation}
\label{V}
U=-V=-{\frac {1}{\sqrt{(x-x_{J})^{2}+(y-y_{J})^{2}+{z}^{2}}}}-(x\ddot x_J+y\ddot y_J)
\end{equation}
is the perturbing gravitational potential. In formula (\ref{V}), the coordinates $(x,y,z)$ of the asteroid can be expressed via Poincar\'e variables by formulas (\ref{coord_transform}), (\ref{P_elements}), and equations of motion of the asteroid in the elliptic orbit
\begin{eqnarray*}
 {{x'}}={{a}}(\cos E-e),\
  {{y'}}={{a}}\sqrt{1-{{e}}^{2}}\,\sin  {{E}} ,\quad
  l=E-e\sin E,
\end{eqnarray*}
where $E$ is the eccentric anomaly of the asteroid. Coordinates $x_J$ and $y_J$ of the planet are prescribed functions of time.  

The double averaged Hamiltonian $\bar F$ is defined as
\begin{equation*}
\bar F =\frac{1}{(2\pi)^2}\int_0^{2\pi}\int_0^{2\pi}F dldl_J.
\end{equation*}
{It is obvious that}
\begin{equation*}
\bar F = -{\frac {\left( 1-\mu \right) ^{2}}{2{L}^{2}}} + \mu {\bar U}= -{\frac {\left( 1-\mu \right) ^{2}}{2{L}^{2}}} - \mu {\bar V},     
\end{equation*}
{where $\mu \bar V=-\mu \bar U$  is the double averaged force function of gravity of the planet:}
\begin{equation*}
\bar U =\frac{1}{(2\pi)^2}\int_0^{2\pi}\int_0^{2\pi}Udldl_J,\quad U={\frac {1}{\sqrt{(x-x_{J})^{2}+(y-y_{J})^{2}+{z}^{2}}}}.
\end{equation*}

By averaging over the mean anomaly, the double averaged Hamiltonian does not depend on $q_1$ anymore, then the canonically conjugate variable $p_1=L$ is the first integral of the double averaged system, which means the first term in $\bar F$ is constant in the double averaged system. Hence, dynamics of variables $p_2$, $q_2$, $p_3$, $q_3$ is described by 2-DOF Hamiltonian system with the Hamiltonian $\mu \bar U$. Introducing ``slow" time $\tau = \mu t$ as an independent variable, then by transformation, the Hamiltonian of the system becomes $\bar{U}$. It is obvious that $\bar U$ depends on the ratio between the semi-major axis of the asteroid and the planet $a$ (we take $a_J=1$) and the eccentricity of the planet $e_J$ as parameters.
 
When the inclination $i=0$, the double averaged spatial ER3BP turns to be the case of planar problem, which corresponds to the invariant plane $p_3=0, q_3=0$. Dynamics in the invariant plane of the system is described by a 1-DOF Hamiltonian system with $p_{2}$ and $q_{2}$. Denote the Hamiltonian of the planar problem by $\bar R(p_2,q_2)$, then $\bar R(p_2,q_2)=\bar U(p_2,q_2,p_3 =0,q_3 =0)$ is independent of $p_3$, $q_3$. According to \citet{akse} and \citet{NSS}, the Hamiltonian of the double averaged planar ER3BP is 
\begin{equation}
\bar R =\frac{1}{(2\pi)^2}\int_0^{2\pi}\int_0^{2\pi}Rdldl_J,\quad R=-{\frac {1}{\sqrt{(x-x_{J})^{2}+(y-y_{J})^{2}}}}.
\end{equation}

We consider the inner case of the double averaged spatial ER3BP. In this case, the distance between the asteroid and the star is much smaller than the distance between the planet and the star. Thus an expansion over the ratio between the semi-major axis of the asteroid and the planet can be established. 

\begin{figure}[htbp]
  \centering
{\includegraphics[width=0.5\textwidth]{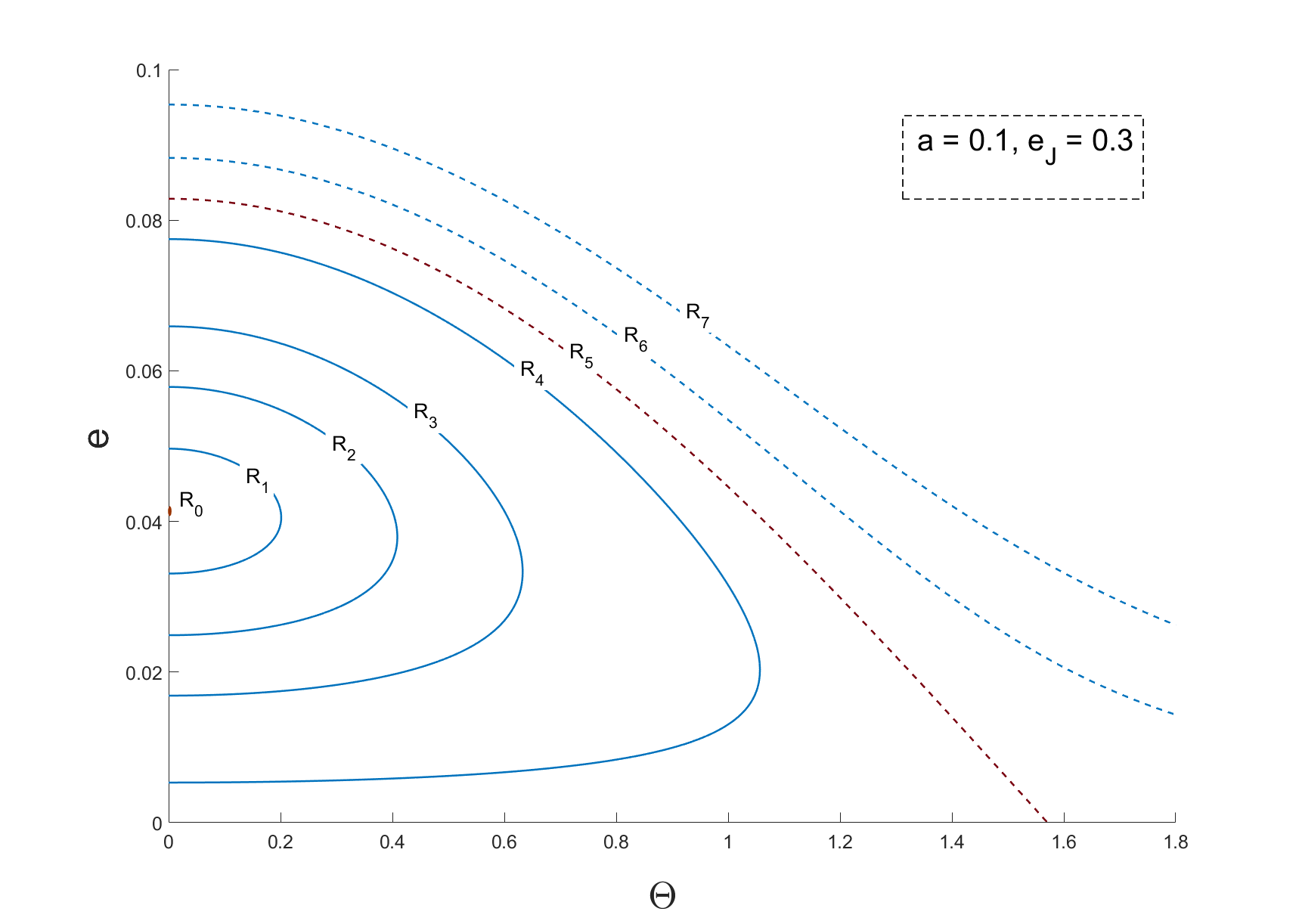}}
\caption{Trajectories of $\Theta$, $e$ in double averaged planar ER3BP.}
\label{fig:Rphase} 
\end{figure}

\begin{figure}[htbp]
  \centering
{\includegraphics[width=0.45\textwidth]{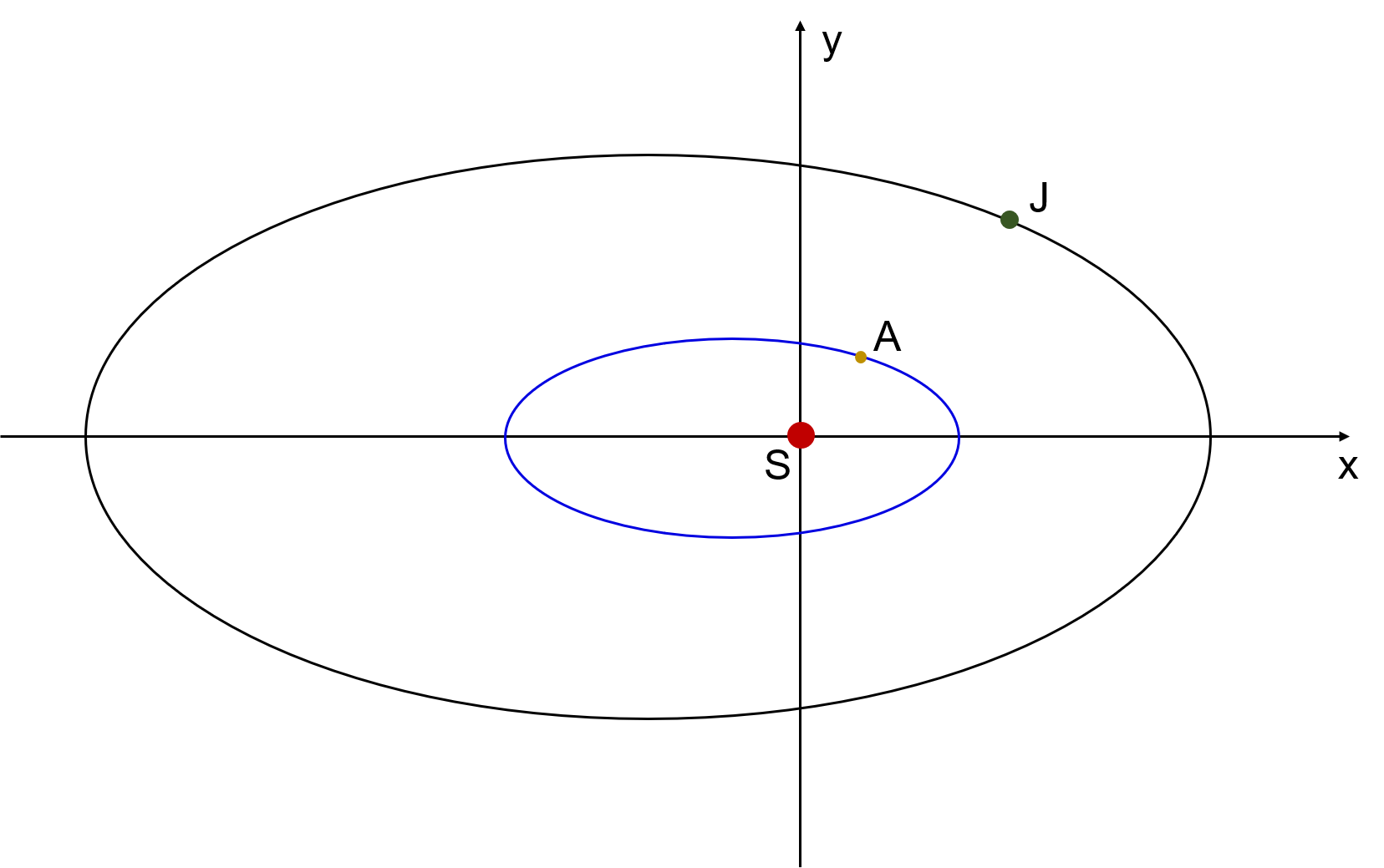}} \quad {\includegraphics[width=0.45\textwidth]{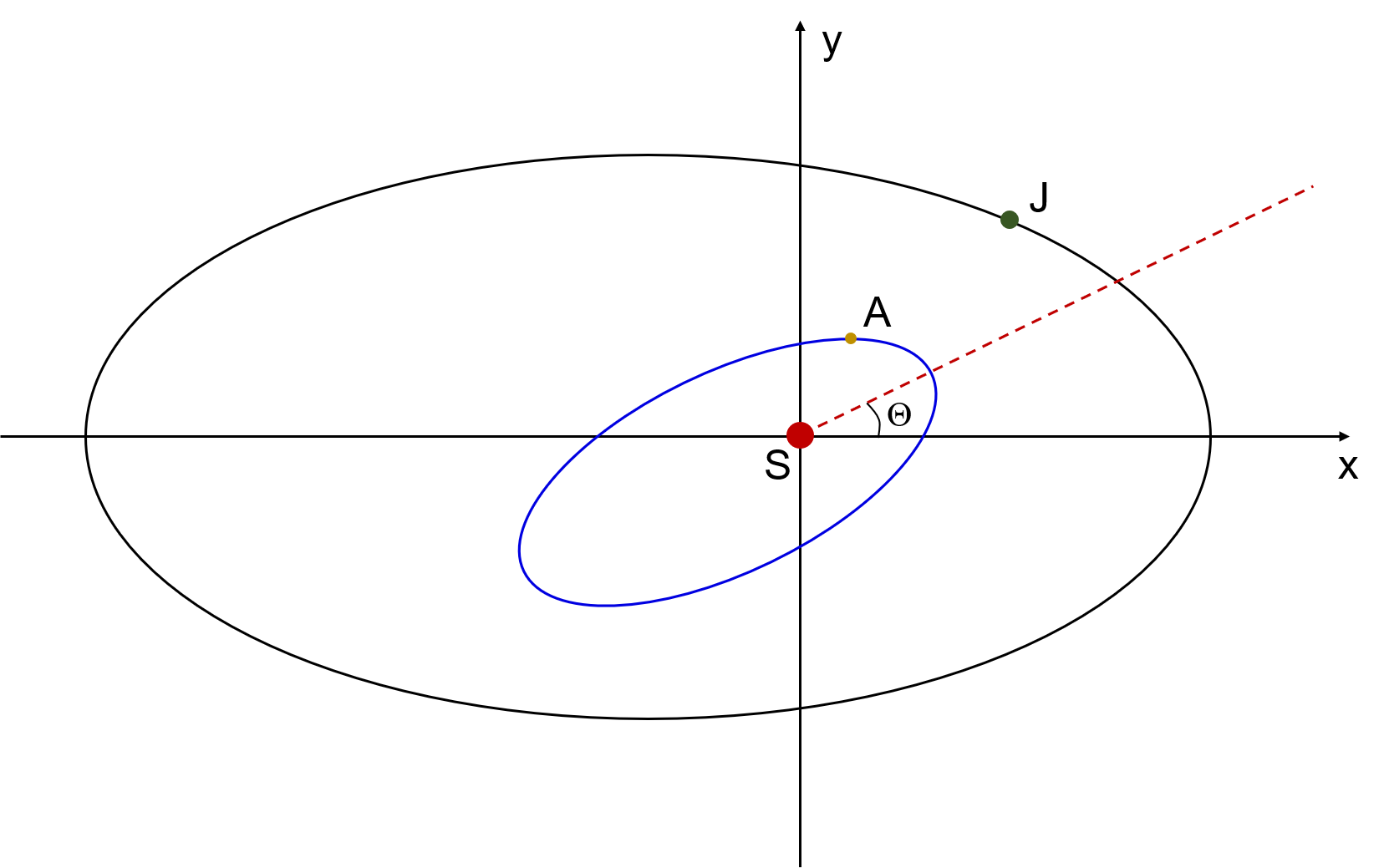}}
\caption{Case of equilibria (apsidal alignment) and periodic orbits.}
\label{fig:alignment} 
\end{figure}

The expression of $\bar R$ truncated in essential orders is 
\begin{equation}
\bar R =-1-\frac{a^{2} \left(3 e^{2}+2\right)}{8 \left(1-e_{J}^{2}\right)^{\frac{3}{2}}}+\frac{15 a^{3} e e_{J} \left(3 e^{2}+4\right) \cos\! \left(\omega +\Omega \right)}{64 \left(1-e_{J}^{2}\right)^{\frac{5}{2}}}.
\end{equation}
Thus, there exists a branch of orbits with respect to $\Theta=\omega+\Omega$ and the eccentricity of the asteroid $e$. Given some values of the eccentricity of the planet $e_{J}$ and the ratio between the semi-major axis of the asteroid and the planet $a$, we obtain figures of the considered orbits ({include equilibria, periodic orbits, separatrix and orbits of circulating}), shown in Fig. \ref{fig:Rphase}, which coincide with the results in \citet{vash_planar}. Assume $a=0.1$ and $e_{J}=0.3$, the equilibria with respect to $R_0\approx 1.002872548$ ($\Theta=0$) is the apsidal alignment case, where $e\approx 0.041367$, see Fig. \ref{fig:alignment} (Left), this case with small inclination $i$ has been studied in \citet{NSS} with a conclusion that the asteroid's orbits are linearly stable. Orbits with types of $R_1=1.002872843$, $R_2=1.002873713$, $R_3=1.002875125$, $R_4=1.002878129$ are periodic orbits, see Fig. \ref{fig:alignment} (Right), while types of $R_6=1.002881952$, $R_7=1.002885001$ are {orbits of circulating}, $R_5\approx 1.002879903$ is the {separatrix of periodic orbits and orbits of circulating}. Periodic orbits only occur in the case of $-\pi/2 < \Theta < \pi/2$.


For spatial problem, denote 

\begin{equation}
M_{1}={\frac {{x}^{2}+{y}^{2}+{z}^{2}}{x_{J}^{2}+y_{J}^{2}}}, \quad  M_{2}=\,-\frac {2\left(xx_{J}+yy_{J}\right)}{x_{J}^{2}+y_{J}^{2}}, \quad M=M_{1}+M_{2}\,. 
\end{equation}
We obtain the approximate formula of the force function of gravity $U$ of the planet for this limiting case. Expansion of  $U$ up to $a^{3}$ has the form 
\begin{equation}
\label{V_expanded}
U=-\frac {1}{\sqrt {x_{J}^{2}+y_{J}^{2}}} \left(1-\frac{1}{2}M+\frac {3}{8} \left(M_{1} M_{2}+ M_{2}^{2}\right)-\frac{5}{16} M_{2}^{3}\right).
\end{equation}

Then we substitute $x$, $x_{J}$, $y$, $y_{J}$ into the above formula and average (\ref{V_expanded}) over $l$, $l_J$ to obtain the corresponding expansion of $\bar{U}$. We consider the apsidal alignment case ($\omega+\Omega=0$) and the general case ($\omega+\Omega=\Theta$) of the problem. The potential function $U$ in the spatial ER3BP is a function of $p_3$, $q_3$ with parameter $\Theta$ and $i$, denoted by $U_{\Theta,i}$. The double averaged potential function $\bar{U}$ is $\bar{U}_{\Theta,i}$. 

{Given spatial perturbations to the equilibria and periodic orbits in planar problem, the Hamiltonian system of 2-DOF evolves to which of 3-DOF, we are interested in the linear stability of the equilibria and the periodic orbits under spatial perturbations.  We consider this kind of spatial problem and study the apsidal alignment case with large inclinations ($U=U_{0,i}$), the general case with small inclinations ($U=U_{\Theta,0^{+}}$), and finally we analyze the general case with large inclinations ($U=U_{\Theta,i}$). }

\section{Stability analysis in apsidal alignment case}

It is established numerically by  \citet{vash_planar}, as well as Fig \ref{fig:Rphase} in this paper, that the double averaged planar restricted elliptic three-body problem  has stationary solutions (equilibria)
\begin{equation}
\label{ss}
q_2=0,\quad p_2=p_{2*},
\end{equation}
for some domains in the plane of parameters $a,e_J$.

The equilibria (\ref{ss}) correspond to the apsidal alignment case in which $\omega+\Omega=0$. We consider a large variation{ (large perturbation)} of the inclination $i$ ($-\pi/2<i<\pi/2$), then the planar problem turns into a spatial problem. It is interesting to discuss the stability of the orbits of the asteroid with respect to an arbitrary inclination $i$.

The expressions of $x$, $y$, $z$ can be simplified as 
\begin{equation}
\label{coord_transform2}
\begin{aligned}
x &=\left( 1-(1-\cos i) \sin ^2 \Omega \right) x'+\left( \cos \Omega \sin \Omega(1-\cos i)  \right) y',
\\ 
 y &=\left( \cos \Omega \sin \Omega(1-\cos i) \right) x' +\left( 1-(1-\cos i) \cos ^2 \Omega  \right) y' ,
\\ 
 z &=\left(\sin  i
  \sin\Omega\right) x'+ \left(\sin i \cos\Omega\right) y'. 
\end{aligned}
\end{equation}

To study the linear stability of the equilibria (\ref{ss}) with respect to arbitrary inclination $i$, we need to consider the quadratic part of $p_3$, $q_3$ in the function $\bar U_{0, i}$ at these equilibria. Substituting $x$, $y$, $z$ into (\ref{V_expanded}) and taking
\begin{equation}
\label{CS1}
\cos (\Omega)=\frac{p_{3}}{\sqrt{G \cdot (1-\cos (i))}},\quad \sin (\Omega)=-\frac{q_{3}}{\sqrt{G \cdot (1-\cos (i))}},
\end{equation}
from (\ref{P_elements}), the function $U=U_{0, i}$ with increasing order of $p_{3}$ and $q_{3}$ is 

\begin{equation}
\label{quadratic1}
U_{0, i}=R_{0}+ W_{0, i} +O(p_3^4+ q_3^4), \quad W_{0, i}= \frac{1}{2}\left(A p_3^2 +2B p_3 q_3+C q_3^2 \right),
\end{equation}
where 

\begin{equation}
\label{ABCR}   
\begin{aligned}
A_{0, i}=&-\frac{2 x' x_{J}+2 \cos \left(i\right) y' y_{J}-\sin\left(i\right)^{2} {y'}^{2}}{\left(x_{J}^{2}+y_{J}^{2}\right)^{\frac{3}{2}} G \left(1-\cos\! \left(i\right)\right)}
\\
&+\frac{3 \left(x'^{2}+y'^{2}\right) \left(x' x_{J}+\cos\left(i\right) y' y_{J}\right)}{\left(x_{J}^{2}+y_{J}^{2}\right)^{\frac{5}{2}} G \left(1-\cos\! \left(i\right)\right)} + O(a^4), \\ 
B_{0, i}=&\frac{\sin\left(i\right)^{2} x' y'-\cos\left(i\right) x' y_{J}-\cos\left(i\right) x_{J} y'+x' y_{J}+x_{J} y'}{\left(x_{J}^{2}+y_{J}^{2}\right)^{\frac{3}{2}} G \left(1-\cos\! \left(i\right)\right)}
\\
&-\frac{3 \left(x'^{2}+y'^{2}\right) \left(x' y_{J}+x_{J} y'\right)}{2 \left(x_{J}^{2}+y_{J}^{2}\right)^{\frac{5}{2}} G} + O(a^4), \\
C_{0, i}=&\frac{x'^{2} \sin\! \left(i\right)^{2}-2 x' x_{J} \cos\! \left(i\right)-2 y' y_{J}}{\left(x_{J}^{2}+y_{J}^{2}\right)^{\frac{3}{2}} G \left(1-\cos\! \left(i\right)\right)}
\\
&+\frac{3 \left(x'^{2}+y'^{2}\right) \left(x' x_{J} \cos\! \left(i\right)+y' y_{J}\right)}{\left(x_{J}^{2}+y_{J}^{2}\right)^{\frac{5}{2}} G \left(1-\cos\! \left(i\right)\right)} + O(a^4),
\end{aligned}
\end{equation}
and $R_{0}$ is zero order terms of $p_{3}$ and $q_{3}$, which is not important here.

We consider the double average value of the coefficients in the quadratic part of  $p_{3}$ and $q_{3}$. Averaging $W_{0, i}$ over the mean anomaly of the asteroid $l$ and the mean anomaly of the planet $l_J$, we have
\begin{equation}
\label{{W}}
\bar{W}_{0, i}=\frac{1}{2}\left(\bar{A}\,{p}_{3}^{2}+2\bar{B} p_3 q_3+\bar{C}\, {q}_{3}^{2}\right),
\end{equation}
where 
\begin{equation}
\begin{aligned}
\label{barAC0}
  & \bar{A}_{0, i}=\frac{1}{{4{\pi }^{2}}}
  \int_{0}^{2\pi }\int_{0}^{2\pi }
  \left( A_{0, i} \right)\text{d}l\text{d}{l}_{J},\\
  & \bar{B}_{0, i}=\frac{1}{{4{\pi }^{2}}}
  \int_{0}^{2\pi }\int_{0}^{2\pi }
  \left( B_{0, i} \right)\text{d}l\text{d}{l}_{J},\\
   & \bar{C}_{0, i}=\frac{1}{{4{\pi }^{2}}}
  \int_{0}^{2\pi }\int_{0}^{2\pi }
  \left(C_{0, i} \right)\text{d}l\text{d}{l}_{J}\\
 \end{aligned}
\end{equation}
are the average values of $A_{0, i}$, $B_{0, i}$ and $C_{0, i}$. 

We have eccentric anomaly $E$ of the asteroid in our formulas of $x$, $y$ and eccentric anomaly $E_{J}$ of the planet in $x_{J}$, $y_{J}$. By Kepler's equation 
\begin{equation}
l=E-e\sin E, \,\,\, l_J=E_J-e\sin E_J,
\end{equation}
then, for any function $f$, we have
\begin{equation}
\label{averageway}
\begin{aligned}
   \frac{1}{{4{\pi }^{2}}} \int_{0}^{2\pi }{}\int_{0}^{2\pi }{f}\text{d}l\text{d}{{l}_{J}}&=\frac{1}{{4{\pi }^{2}}} \int_{0}^{2\pi }{}\int_{0}^{2\pi }{f}        
   \frac{\text{d}{{l}}}{\text{d}{{E}}}
   \frac{\text{d}{{l}_{J}}}{\text{d}{{E}_{J}}}     {\text{d}E}    \text{d}{{E}_{J}} \\ 
 & =\frac{1}{{4{\pi }^{2}}} \int_{0}^{2\pi }{}\int_{0}^{2\pi }{f}\left( 1-e\cos  E  \right)\left( 1-{{e}_{J}}\cos  {{E}_{J}}  \right)\text{d}E\text{d}{{E}_{J}}. 
\end{aligned}
\end{equation}  
Thus the averaging values $\bar{A}_{0, i}$, $\bar{B}_{0, i}$ and $\bar{C}_{0, i}$ can be written as 

\begin{equation}
\begin{aligned}
\label{barAC}
  & \bar{A}_{0, i}=\frac{1}{{4{\pi }^{2}}}
  \int_{0}^{2\pi }\int_{0}^{2\pi }
  \left(A_{0, i} \right)\left( 1-e\cos  E  \right)\left( 1-{{e}_{J}}\cos {{E}_{J}}  \right)\text{d}E\text{d}{E}_{J},\\
  & \bar{B}_{0, i}=\frac{1}{{4{\pi }^{2}}}
  \int_{0}^{2\pi }\int_{0}^{2\pi }
  \left( B_{0, i} \right)\left( 1-e\cos  E  \right)\left( 1-{{e}_{J}}\cos {{E}_{J}}  \right)\text{d}E\text{d}{E}_{J},\\
   & \bar{C}_{0, i}=\frac{1}{{4{\pi }^{2}}}
  \int_{0}^{2\pi }\int_{0}^{2\pi }
  \left( C_{0, i} \right)\left( 1-e\cos  E  \right)\left( 1-{{e}_{J}}\cos {{E}_{J}}  \right)\text{d}E\text{d}{E}_{J}.
 \end{aligned}
\end{equation}
We substitute $x'$, $y'$, $x_{J}$, $y_{J}$ into the above formula. It can be calculated that the value of $\bar{B}_{0, i}$ is always $0$, and 

\begin{equation}
\begin{aligned}
\label{barAC1}
& \bar{A}_{0, i}=\frac{\left(1-e^{2}\right) \sin\! \left(i\right)^{2} a^{2}}{4 G \left(1-\cos\! \left(i\right)\right) \left(1-e_{J}^{2}\right)^{\frac{3}{2}}}-\frac{15 e_{J} \left(4+3 e^{2}\right) e a^{3}}{16 G \left(1-\cos\! \left(i\right)\right) \left(1-e_{J}^{2}\right)^{\frac{5}{2}}},\\
& \bar{C}_{0, i}=\frac{\left(1+4 e^{2}\right) \sin\! \left(i\right)^{2} a^{2}}{4 G \left(1-\cos\! \left(i\right)\right) \left(1-e_{J}^{2}\right)^{\frac{3}{2}}}-\frac{15 e_{J} \cos\! \left(i\right) \left(4+3 e^{2}\right) e a^{3}}{16 G \left(1-\cos\! \left(i\right)\right) \left(1-e_{J}^{2}\right)^{\frac{5}{2}}}.\\
 \end{aligned}
\end{equation}

\begin{figure}[htbp]
  \centering
{\includegraphics[width=0.45\textwidth]{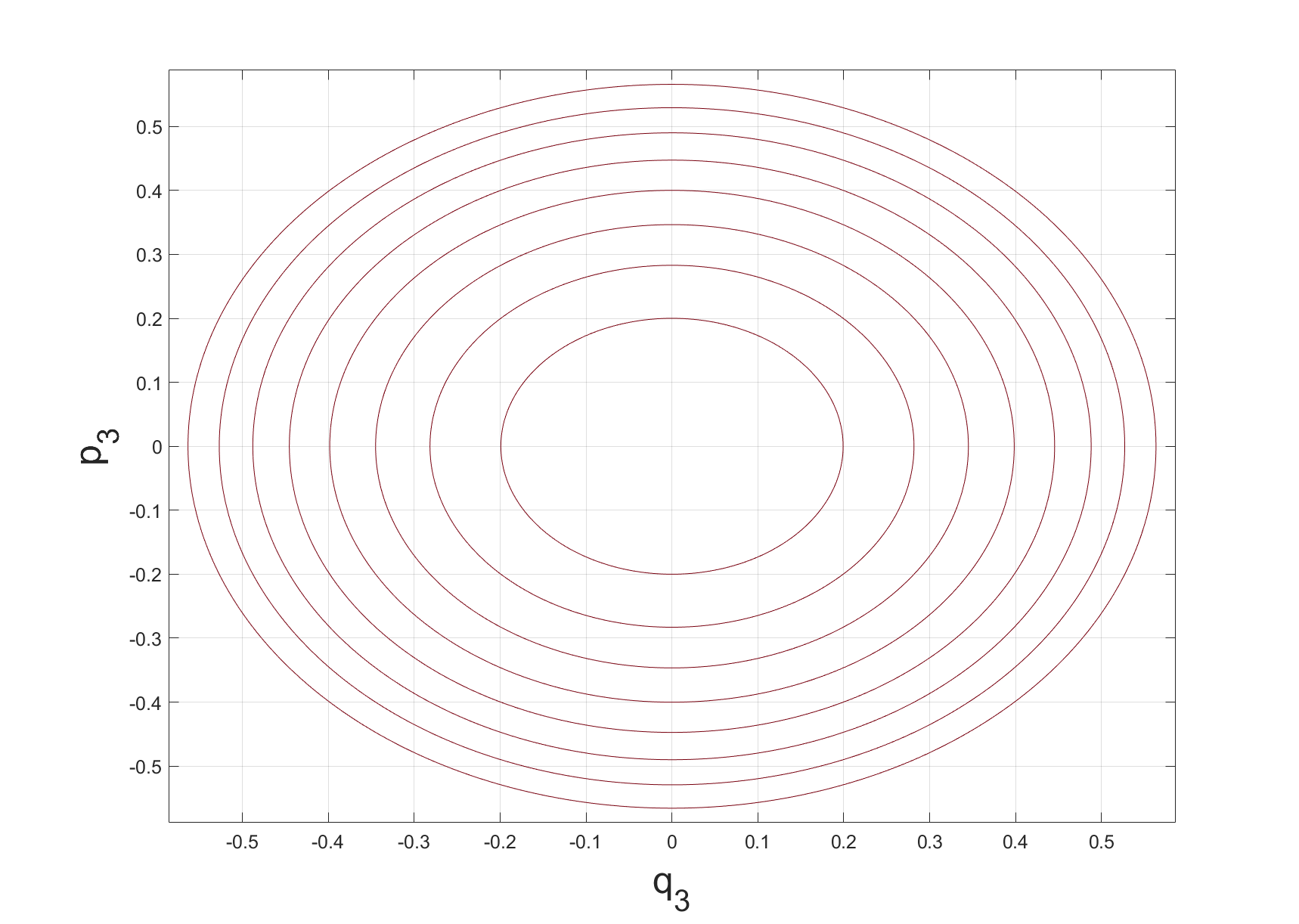}}
\caption{Phase portrait of $p_3$, $q_3$ in apsidal alignment case.}
\label{fig:Phase0} 
\end{figure}

{In mathematics, when $a$ is small, bounded properties of coefficients of terms of $a^3$ is required to guarantee that terms of $a^3$ are much smaller than terms of $a^2$. For $e_{J} \nrightarrow 1$ and $i \nrightarrow 0$, coefficients of $a^3$ in (\ref{barAC1}) are bounded. Then it is obvious that $\bar{A}_{0, i}>0$ and $\bar{C}_{0, i}>0$ for sufficiently small $a$, thus $\bar W_{0, i}$ is a positive definite quadratic form.} Hence, stable equilibria of the double averaged planar restricted elliptic three-body problem are stable in the linear approximation as equilibria of the double averaged spatial restricted elliptic three-body problem for all values of parameters. Numerical simulation is shown in Fig. \ref{fig:Phase0} with $e_J=0.3$, $a=0.1$, and the corresponding equilibrium is $\Theta=0$, $e=0.041367$.

It is known that the double averaged planar ER3BP corresponds to the invariant plane $i=0$. However, the asteroid can not cross the invariant plane without external forces {in the quadrupole model}. This is because that the vertical angular momentum $H$ is also conserved in the system in addition to the total energy \citep{naoz2013}. Since $H=G \cos i$ is a constant, and $G$ is always positive, the test particle should always be in a prograde orbit or always be in a retrograde orbit so that $\cos i$ does not change its sign. Thus, given a large perturbation of the inclination $i$ in the above problem, the orbits of the asteroid will return to the invariant plane with infinite time. Although the system can not reach the equilibrium in finite time, a situation of small inclination $i$ is able to appear, which is the case has been studied in \citet{NSS} with a result of linear stability.

\section{Stability analysis in small perturbation case}

The equilibria (\ref{ss}) correspond to the apsidal alignment case in which $\omega+\Omega=0$. \citet{vash_planar} and Fig \ref{fig:Rphase} in this paper demonstrated that there are period orbits and other orbits around the equilibria (\ref{ss}), which correspond to $\Theta=\omega+\Omega\neq 0$. Thus, consideration of general cases (i.e. $\Theta \neq 0$) is useful to establish the linear stability of the asteroid's orbits. 

We consider the small perturbation case initially, which means the inclination $i$ is regarded as the small perturbation. In this case, the formula of $x$, $y$ and $z$ can be obtained with $i=0$ and $\omega+\Omega=\Theta$, thus we have $z=0$ and  
\begin{equation}
\label{coord_transform2'}
\begin{aligned}
x &=\cos(\Theta) x'-\sin(\Theta) y',
\\ 
y &=\sin(\Theta) x'+\cos(\Theta) y'.
 \end{aligned}
\end{equation}

Formula (\ref{CS1}) can be expanded as 
\begin{equation}
\cos (\Omega)=\frac{p_{3}}{\sqrt{G \cdot i^2/2}},\quad \sin (\Omega)=-\frac{q_{3}}{\sqrt{G \cdot i^2/2}}.
\end{equation}
Substituting $x$, $y$ and $z=0$ into (\ref{V_expanded}), the function $U=U_{\Theta,0^{+}}$ with increasing order of $p_{3}$ and $q_{3}$ can be written as
\begin{equation}
\label{quadratic2}
U_{\Theta, 0^{+}}=R_{\Theta}+ W_{\Theta, 0^{+}} +O(p_3^4+ q_3^4), \quad W_{\Theta, 0^{+}}= \frac{1}{2}\left(A_{\Theta, 0^{+}} p_3^2 +2B_{\Theta, 0^{+}} p_3 q_3+C_{\Theta, 0^{+}} q_3^2 \right),
\end{equation}
where 
\begin{equation}
\label{ABCR0}   
\begin{aligned}
A_{\Theta, 0^{+}}&=\,\frac{y{{y}_{J}}}{G {\left[ {\left( x-{{x}_{J}} \right)}^{2}+{{\left( y-{{y}_{J}} \right)}^{2}} \right]}^{{3}/{2}}}, 
 \\ 
B_{\Theta, 0^{+}}&=\,\frac{\left( x{{y}_{J}}+y{{x}_{J}} \right)}{2G {{\left[ {{\left( x-{{x}_{J}} \right)}^{2}}+{{\left( y-{{y}_{J}} \right)}^{2}} \right]}^{{3}/{2}}}}, 
\\
C_{\Theta, 0^{+}}&=\, \frac{x{{x}_{J}}}{ G {\left[ {{\left( x-{{x}_{J}} \right)}^{2}}+{{\left( y-{{y}_{J}} \right)}^{2}} \right]^{{3}/{2}}}}.
\end{aligned}
\end{equation}
Calculations are similar to the procedure in \citet{NSS}.

We consider the double average value of the coefficients in the quadratic form of  $p_{3}$ and $q_{3}$, and average $W_{\Theta, 0^{+}}$ with the same procedure in (\ref{averageway}), then 
\begin{equation}
\label{{W}_{in} 0}
\bar{W}_{\Theta, 0^{+}}=\frac{1}{2}\left(\bar{A}_{\Theta, 0^{+}}\,{{p}_{3}}^{2}+2\bar{B}_{\Theta, 0^{+}} p_3 q_3+\bar{C}_{\Theta, 0^{+}}\, {{q}_{3}}^{2}\right),
\end{equation}
where $\bar{A}_{\Theta, 0^{+}}$, $\bar{B}_{\Theta, 0^{+}}$ and $\bar{C}_{\Theta, 0^{+}}$ are the average values of $A_{\Theta, 0^{+}}$, $B_{\Theta, 0^{+}}$ and $C_{\Theta, 0^{+}}$ respectively by the averaging method in formula (\ref{barAC0}) and (\ref{barAC}). 

Substituting $x'$, $y'$, $x_{J}$, $y_{J}$ into the double average value of the coefficients, we obtain 

\begin{equation}
\label{ABCtheta0}   
\begin{aligned}
\bar{A}_{\Theta, 0^{+}}=&\frac{3 \left(1+4 e^{2}-5 e^{2} \cos \left(\Theta \right)^{2} \right)}{4 G \left(1-e_{J}^{2}\right)^{\frac{3}{2}}}\,a^{2} 
\\
&+\frac{15\, e\, e_{J} \cos \left(\Theta \right) \left(70 \cos \! \left(\Theta \right)^{2} e^{2}-60 e^{2}-10\right)}{64 G \left(1-e_{J}^{2}\right)^{\frac{3}{2}}}\,a^{3} + O(a^4), 
\\ 
\bar{B}_{\Theta, 0^{+}}=&\frac{15 e^{2} \cos\! \left(\Theta \right) \sin\! \left(\Theta \right)}{4 G \left(-e_{J}^{2}+1\right)^{\frac{3}{2}}}\,a^{2}
\\
&-\frac{15 e \sin\! \left(\Theta \right) e_{J} \left(140 \cos\! \left(\Theta \right)^{2} e^{2}-17 e^{2}+24\right)}{128 G \left(1-e_{J}^{2}\right)^{\frac{3}{2}}}\,a^{3} + O(a^4), 
\\ 
\bar{C}_{\Theta, 0^{+}}=&\frac{3 \left(1- e^{2}+5 e^{2} \cos \left(\Theta \right)^{2} \right)}{4 G \left(1-e_{J}^{2}\right)^{\frac{3}{2}}}\,a^{2}
\\
&-\frac{15\, e\, e_{J} \cos \left(\Theta \right) \left(70 \cos\! \left(\Theta \right)^{2} e^{2}-27 e^{2}+34\right)}{64 G \left(1-e_{J}^{2}\right)^{\frac{3}{2}}}\,a^{3} + O(a^4). 
\end{aligned}
\end{equation}

{In above formulas, as $a$ is small enough, for $e_{J} \nrightarrow 1$, coefficients of $a^3$ in (\ref{ABCtheta0}) are bounded, thus terms of higher order of $a$ are omitted.} It is obvious that $\bar{A}_{\Theta, 0^{+}}>0$ and $\bar{C}_{\Theta, 0^{+}}>0$. The determinant of the quadratic form (\ref{{W}_{in} 0}) is
\begin{equation}
\label{DT0}
\begin{aligned}
&D_{\Theta, 0^{+}}= \det\left[\begin{array}{cc}
\bar{A}_{\Theta, 0^{+}} & \bar{B}_{\Theta, 0^{+}} 
\\
\bar{B}_{\Theta, 0^{+}} & \bar{C}_{\Theta, 0^{+}} 
\end{array}\right] 
\\
=&\frac{9 \left(1+3 e^{2}-4 e^{4}\right)}{16\, G^{2} \left(1-e_{J}^{2}\right)^{3}}\,a^{4}+\frac{45 e_{J} \left(83 e^{4}-39 e^{2}-44\right) e \cos\! \left(\Theta \right)}{256\, G^{2} \left(1-e_{J}^{2}\right)^{4}}\,a^{5}
\\
&+\frac{225 \left(\begin{aligned}1431 \cos\! \left(\Theta \right)^{2} e^{4}+456 \cos\! \left(\Theta \right)^{2} e^{2}+289 e^{4}\\-1936 \cos\! \left(\Theta \right)^{2}-816 e^{2}+576 \end{aligned}\right) e^{2} e_{J}^{2}}{16384\, G^{2} \left(1-e_{J}^{2}\right)^{5}}\,a^{6}.
\end{aligned}
\end{equation}

{Similarly, when $a$ is small, coefficients of $a^5$ and $a^6$ in (\ref{DT0}) are bounded as far as $e_{J} \nrightarrow 1$.} For small values of $a$, the sequential principal minor $\bar{A}_{\Theta, 0^{+}}>0$ and $D_{\Theta, 0^{+}} >0$, thus (\ref{{W}_{in} 0}) is a positive definite quadratic form. These conditions provide linear stability of the orbits of the asteroid in small perturbation case, which means stable periodic orbits of the double averaged planar restricted elliptic three-body problem are stable in the linear approximation as periodic orbits of the double averaged spatial restricted elliptic three-body problem with small inclination for all values of parameters.

\section{Expansion of the problem for general inner case}

The spatial perturbation of the asteroid's orbits may not be small, as we discussed in previous section, large variations {(large perturbation)} of the inclination $i$ are of great interest. In this section, we consider the problem with large perturbation, i.e. inclination $i$ is not small. 

For general situation, we have $\omega+\Omega=\Theta$ in the formula of $x$,$y$ and $z$. Thus
\begin{equation}
\label{coord_transform3}
\begin{aligned}
x =&\left( \cos\Omega\cos(\Theta-\omega) -\cos  i\sin\Omega  \sin(\Theta-\omega)  \right) x'
\\
&+\left( -\cos\Omega\sin  (\Theta-\omega)  -\cos i\sin\Omega \cos(\Theta-\omega)  \right) y',
\\ 
 y =&\left( \sin \Omega  \cos(\Theta-\omega) +\cos i\cos\Omega \sin(\Theta-\omega)  \right) x' 
 \\
 &+\left( -\sin\Omega\sin(\Theta-\omega) +\cos  i\cos\Omega  \cos(\Theta-\omega)  \right) y' ,
\\ 
 z =&\left(\sin  i  \sin(\Theta-\omega)\right) x'+\left(\sin i \cos(\Theta-\omega)\right) y'. 
 \end{aligned}
\end{equation}

We still take 
$$
\cos (\Omega)=\frac{p_{3}}{\sqrt{G \cdot (1-\cos (i))}},\quad \sin (\Omega)=-\frac{q_{3}}{\sqrt{G \cdot (1-\cos (i))}}
$$
from (\ref{P_elements}). Substituting $x$, $y$, $z$ into (\ref{V_expanded}), the function $U_{\Theta}$ with increasing order of $p_{3}$ and $q_{3}$ can be written as
\begin{equation}
\label{quadratic3}
U_{\Theta,i}=R_{\Theta}+ W_{\Theta,i} +O(p_3^4+ q_3^4), \quad W_{\Theta,i}= \frac{1}{2}\left(A_{\Theta,i} p_3^2 +2B_{\Theta,i} p_3 q_3+C_{\Theta,i} q_3^2 \right),
\end{equation}
where $A_{\Theta,i}$, $B_{\Theta,i}$ and $C_{\Theta,i}$ are coefficients of $p_3^2$, $p_3 q_3$ and $q_3^2$ respectively, which can be calculated with a similar procedure as (\ref{ABCR}).




We consider the double average value of the coefficients in the quadratic form of  $p_{3}$ and $q_{3}$, and average $W_{\Theta,i}$ with the same procedure in (\ref{averageway}), then 
\begin{equation}
\label{{W}_{in}}
\bar{W}_{\Theta,i}=\frac{1}{2}\left(\bar{A}_{\Theta,i}\,{{p}_{3}}^{2}+2\bar{B}_{\Theta,i} p_3 q_3+\bar{C}_{\Theta,i}\, {{q}_{3}}^{2}\right),
\end{equation}
where 

\begin{equation}
\begin{aligned}
\label{bartheta}
  & \bar{A}_{\Theta,i}=\frac{1}{{4{\pi }^{2}}}
  \int_{0}^{2\pi }\int_{0}^{2\pi }
  \left(A_{\Theta,i} \right)\left( 1-e\cos  E  \right)\left( 1-{{e}_{J}}\cos {{E}_{J}}  \right)\text{d}E\text{d}{E}_{J},\\
  & \bar{B}_{\Theta,i}=\frac{1}{{4{\pi }^{2}}}
  \int_{0}^{2\pi }\int_{0}^{2\pi }
  \left( B_{\Theta,i} \right)\left( 1-e\cos  E  \right)\left( 1-{{e}_{J}}\cos {{E}_{J}}  \right)\text{d}E\text{d}{E}_{J},\\
   & \bar{C}_{\Theta,i}=\frac{1}{{4{\pi }^{2}}}
  \int_{0}^{2\pi }\int_{0}^{2\pi }
  \left( C_{\Theta,i} \right)\left( 1-e\cos  E  \right)\left( 1-{{e}_{J}}\cos {{E}_{J}}  \right)\text{d}E\text{d}{E}_{J}\\
 \end{aligned}
\end{equation}
are the average values of $A_{\Theta,i}$, $B_{\Theta,i}$ and $C_{\Theta,i}$.

Substituting $x'$, $y'$, $x_{J}$, $y_{J}$ into the double average value of the coefficients, we obtain 

\begin{equation}
\label{ABCtheta}   
\begin{aligned}
\bar{A}_{\Theta,i}&=\frac{\sin \! \left(i \right)^{2}\left(1+4 e^{2}-5 \cos \! \left(\Theta \right)^{2} e^{2}\right) }{4 G \left(1-e_{J}^{2}\right)^{\frac{3}{2}} \left(1-\cos\left(i \right)\right)}\, a^{2} -\frac{15 \cos\! \left(\Theta \right)  e_{J} e  \left(3 e^{2}+4\right) }{16 G \left(1-\cos \left(i \right)\right) \left(1-e_{J}^{2}\right)^{\frac{5}{2}}}\, a^{3} , \\ 
\bar{B}_{\Theta,i}&=\frac{5 \sin\! \left(i\right)^{2} \cos\! \left(\Theta \right) \sin\! \left(\Theta \right) e^{2}}{4 G \left(1-e_{J}^{2}\right)^{\frac{3}{2}} \left(1-\cos\left(i\right)\right)}\, a^{2}+\frac{15 \sin\! \left(\Theta \right) e_{J} e \left(3 e^{2}+4\right)}{32 G \left(1-e_{J}^{2}\right)^{\frac{5}{2}}}\, a^{3}, \\ 
\bar{C}_{\Theta,i}&=\frac{\sin\! \left(i\right)^{2}  \left(5 \cos\! \left(\Theta \right)^{2} e^{2}-e^{2}+1\right)}{4 G \left(1-e_{J}^{2}\right)^{\frac{3}{2}} \left(1-\cos \left(i\right)\right)}\, a^{2}-\frac{15 \cos \left(i \right) \cos\! \left(\Theta \right)  e_{J} e   \left(3 e^{2}+4\right) }{16 G \left(1-\cos \left(i\right)\right) \left(1-e_{J}^{2}\right)^{\frac{5}{2}}}\, a^{3}. 
\end{aligned}
\end{equation}
The determinant of the quadratic form (\ref{{W}_{in}}) is
\begin{equation}
\label{D} 
\begin{aligned}
D_{\Theta,i}=& \det\left[\begin{array}{cc}
\bar{A}_{\Theta,i} & \bar{B}_{\Theta,i} 
\\
 \bar{B}_{\Theta,i} & \bar{C}_{\Theta,i} 
\end{array}\right]=\frac{\left(1+3 e^{2}-4 e^{4}\right) \left(1+\cos \left(i\right)\right)^{2}}{16 \left(1-e_{J}^{2}\right)^{3} G^{2}}\,a^{4} 
\\
&+\frac{15 \sin\left(i \right)^{2} \cos \left(\Theta \right) \left(-4 e^{2}-1+\cos \left(i \right) e^{2}-\cos\! \left(i \right)\right) e_{J} \left(3 e^{2}+4\right) e}{64 G^{2} \left(1-\cos \left(i \right)\right)^{2} \left(1-e_{J}^{2}\right)^{4}}\,a^{5}
\\
&-\frac{225 \left(\begin{aligned} -\cos \left(i \right)^{2} \sin \left(\Theta \right)^{2}+2 \cos \left(i \right) \cos \left(\Theta \right)^{2}\\ -\sin \left(\Theta \right)^{2}+2 \cos \left(i \right)\end{aligned} \right) e_{J}^{2} \left(3 e^{2}+4\right)^{2} e^{2}}{1024 G^{2} \left(1-\cos \left(i \right)\right)^{2} \left(1-e_{J}^{2}\right)^{5}}\,a^{6}.
\end{aligned}
\end{equation}

For $e_{J} \nrightarrow 1$ and $i \nrightarrow 0$, coefficients of $a^3$ in (\ref{ABCtheta}) and coefficients of $a^5$ and $a^6$ in (\ref{D}) are bounded, thus for sufficient small $a$, higher order terms of $a$ can be omitted. In this case, for any values of $e$, $\bar{A}_{\Theta,i}>0$,  $\bar{C}_{\Theta,i}>0$ and $D_{\Theta,i}>0$ hold. The sequential principal minor $\bar{A}_{\Theta,i}>0$ and $D_{\Theta,i}=\bar{A}_{\Theta,i} \bar{C}_{\Theta,i}-\bar{B}_{\Theta,i}^2 >0$, thus (\ref{{W}_{in}}) is a positive definite quadratic form, which guaranteed the linear stability of the orbits of the asteroid. Hence, given a stable periodic orbit of the double averaged planar restricted elliptic three-body problem, linear stability is kept for the periodic orbits of the double averaged spatial restricted elliptic three-body problem with large perturbation over inclination $i$ for all values of parameters. 

The major term of $D_{\Theta,i}$ (i.e. term of $a^4$) is in absence of variable $\Theta$. Thus, the positive definiteness of the determinant does not depend on $\Theta=\omega+\Omega$ in inner case. When $\Theta=0$, the problem corresponds to the apsidal alignment case, and the conclusion coincides with the results in the previous section.

\begin{figure}[htbp]
  \centering
 {\includegraphics[width=0.45\textwidth]{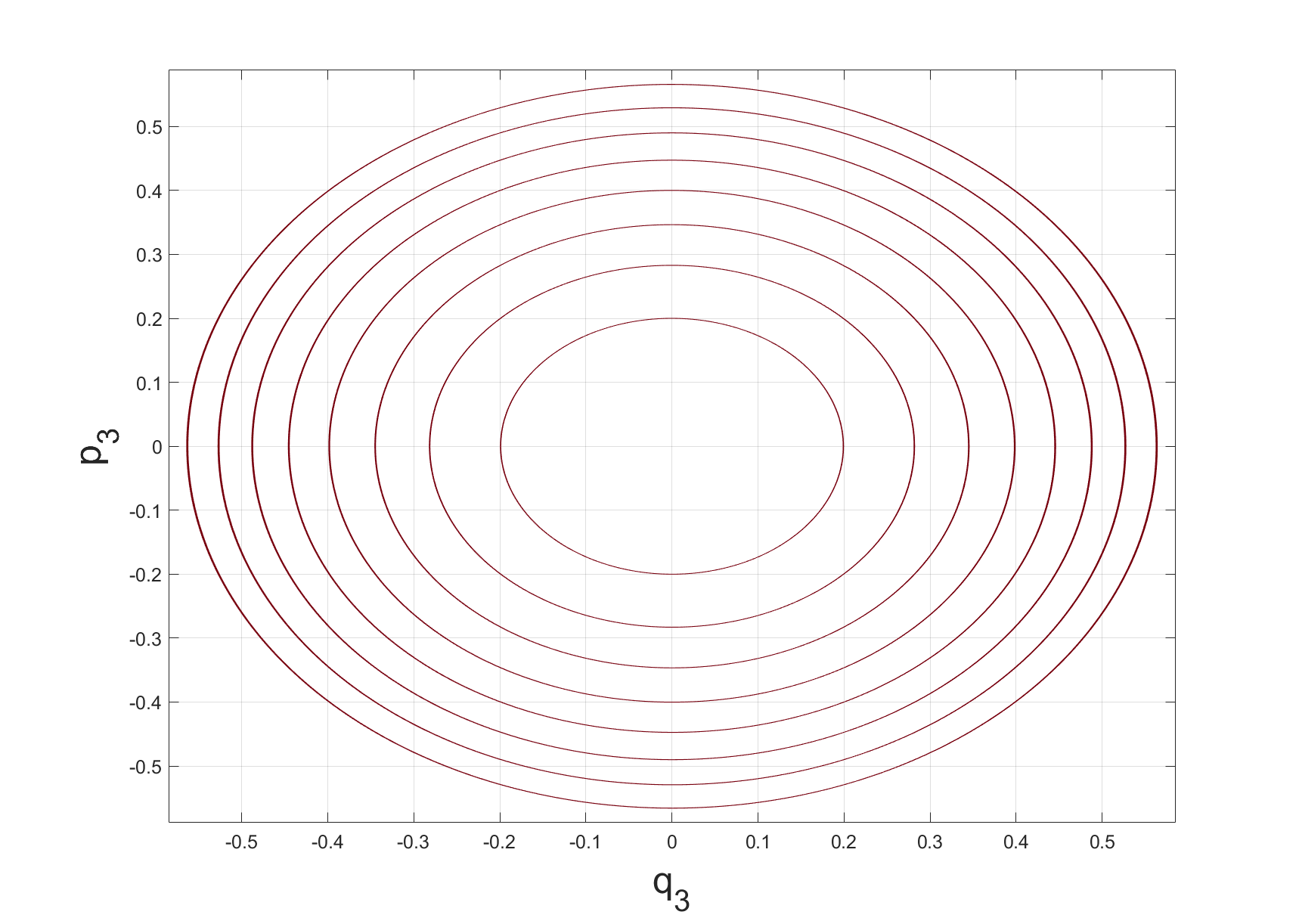}}\quad
{\includegraphics[width=0.45\textwidth]{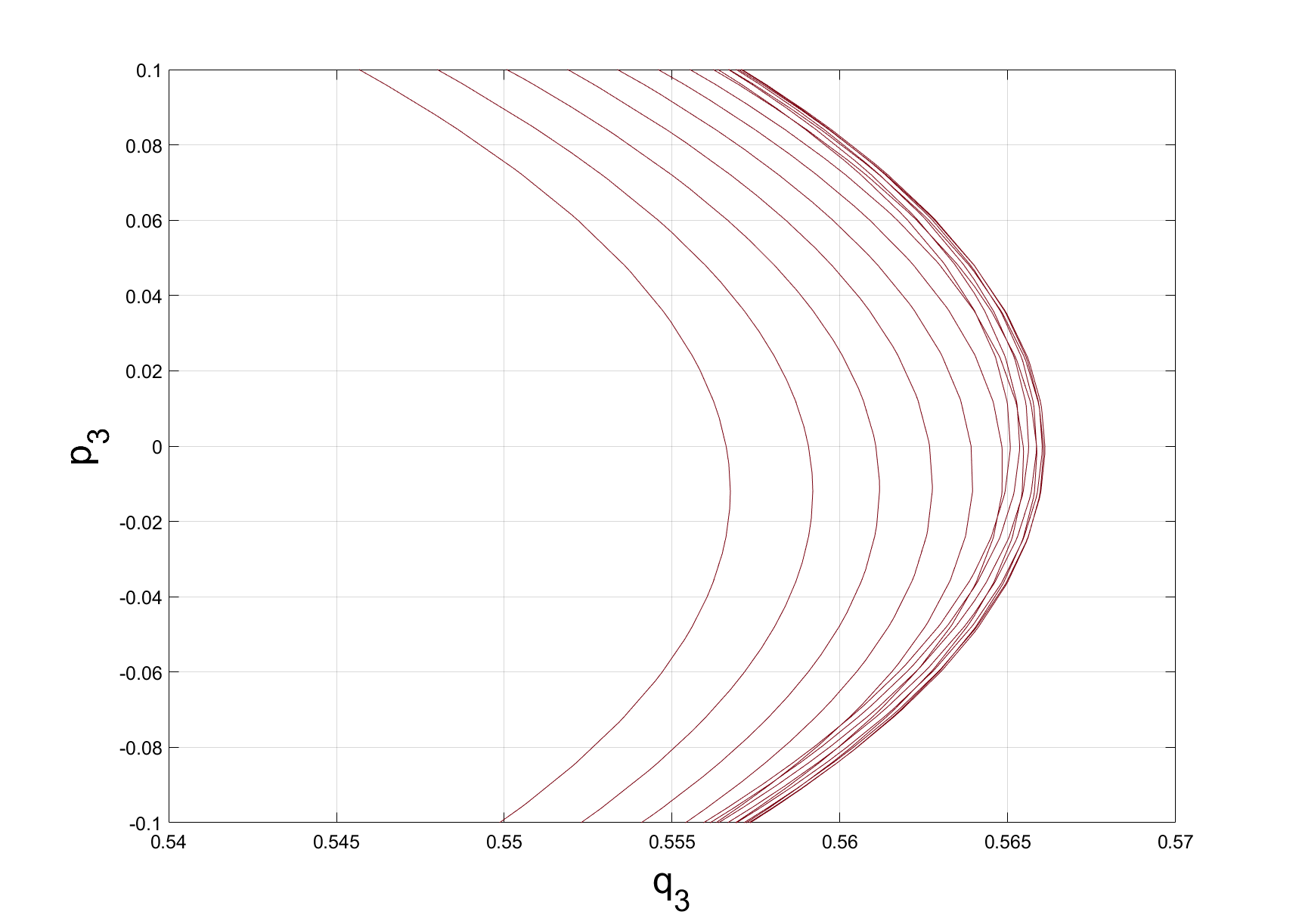}}
\caption{Phase portrait of $p_3$, $q_3$ with periodic orbits of $\Theta$ and $e$.}
\label{fig:Phase12} 
\end{figure}

\begin{figure}[htbp]
  \centering
{\includegraphics[width=0.45\textwidth]{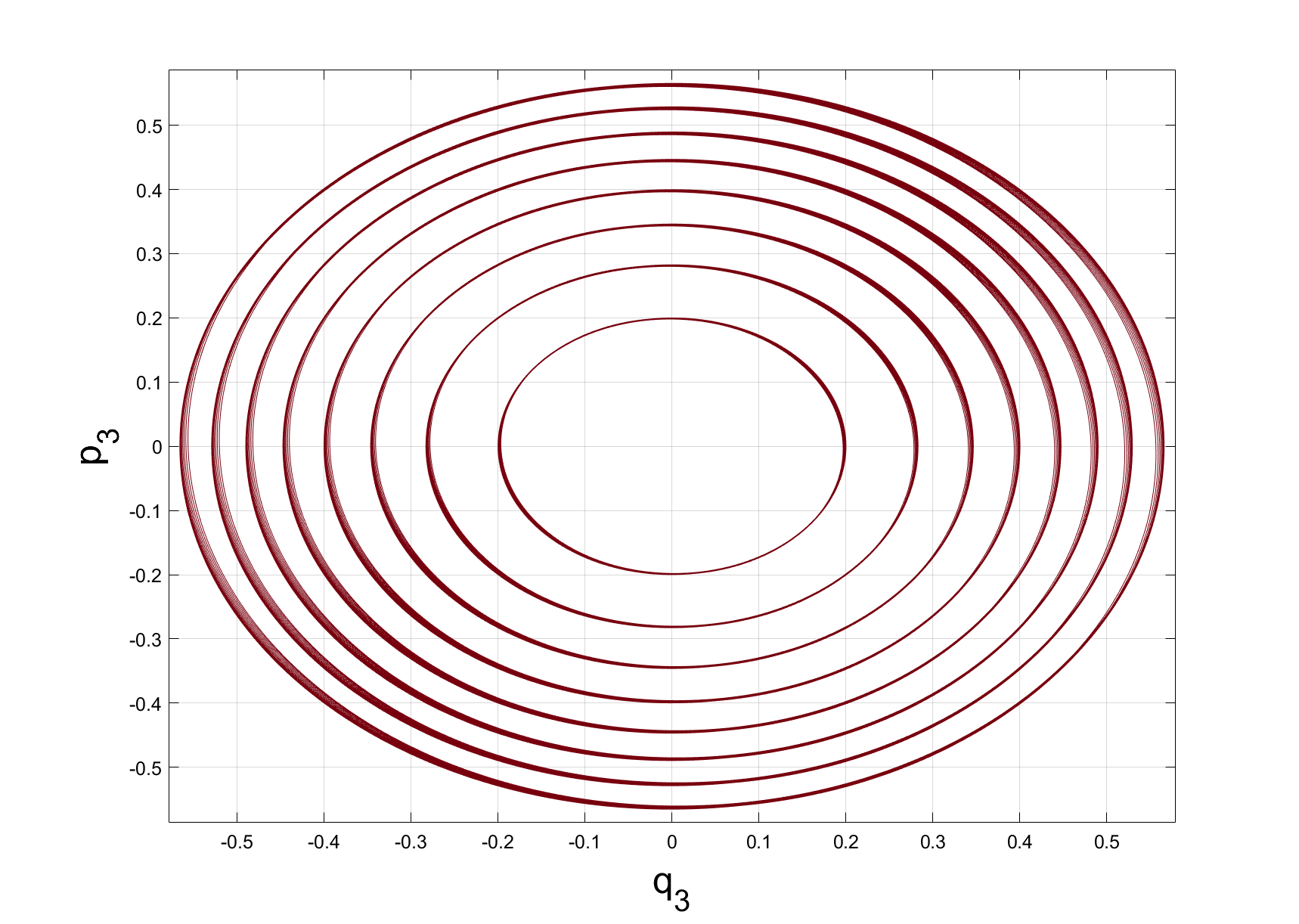}}\quad
{\includegraphics[width=0.45\textwidth]{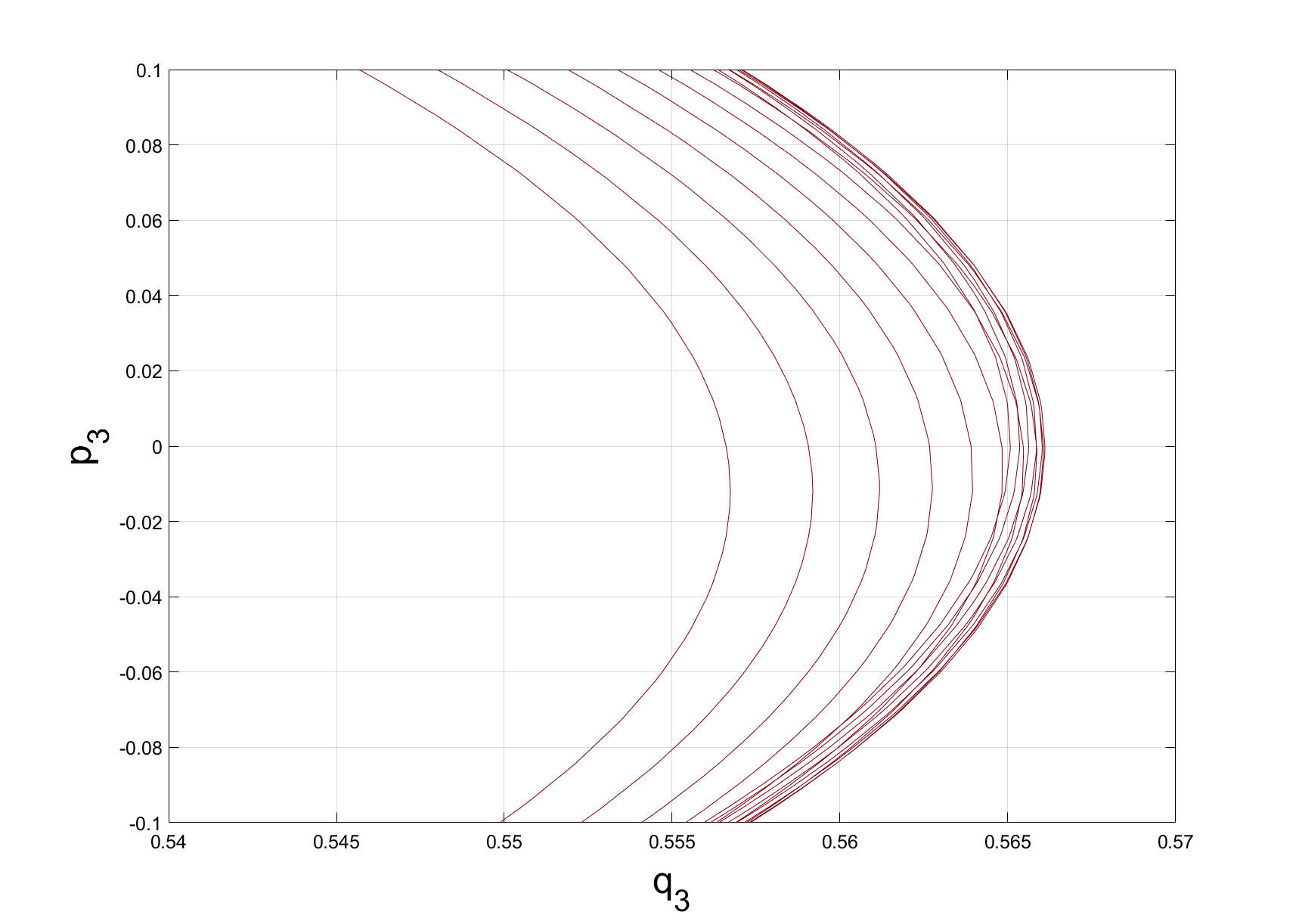}}
\caption{Phase portrait of $p_3$, $q_3$ with other orbits of $\Theta$ and $e$.}
\label{fig:Phase34} 
\end{figure}

Numerical examples of phase portraits of the system for canonical variables $p_3$ and $q_3$ with periodic orbits and other orbits of parameters $\Theta$ and $e$ are shown in Fig. \ref{fig:Phase12} and Fig. \ref{fig:Phase34} respectively with $e_J=0.3$ and $a=0.1$. Fig. \ref{fig:Phase12} (Left) is a description of stability of the periodic orbit corresponding to $R_1=1.002872843$, Fig. \ref{fig:Phase34} (Left) shows the stability of orbit corresponding to $R_7=1.002885001$. Fig. \ref{fig:Phase12} (Right) and Fig. \ref{fig:Phase34} (Right) are the enlarged fragments of the graphs presented on the left panels, which demonstrate that the changing of values of $e$ in the same orbit cause small influence of the phase curves, hence it shows clearly in right panels that the linear stability of considered orbits are guaranteed. 

Similar to the explanation in the apsidal alignment case, the orbits of the asteroid can not cross the invariant plane and will return to the invariant plane in finite time. Thus, in finite time, the problem will degenerate to the small perturbation case, which is studied in the previous section.

We only considered the double averaged spatial restricted elliptic three-body problem with the inclination $-\pi/2<i<\pi/2$, for values of $i$ between $(\pi/2,\pi)$ and $(-\pi/2,-\pi/2)$, one can choose the corresponding equilibria in the apsidal contrary case, and perform a similar calculation. 

We only consider the linear stability in this work, Lidov-Kozai effect resulted in large coupled periodic changes of the eccentricity and the inclination of a natural or artificial celestial body. Thus, the asteroid's orbits may not prevent the stability if there are nonlinear instability or resonances instability between frequencies of the eccentricity and the inclination. Most natural or artificial celestial bodies have small eccentricities and inclinations, the mentioned resonances instability hardly occur. Thus, the study of linear stability of periodic orbits around the apsidal alignment equilibria  (\ref{ss}) is of great importance. Further discussion is needed to determine the nonlinear stability of the problem, which is an open question. The analysis of the nonlinear stability as well as the resonance cases can be a direction of future work.

\section{Conclusion}

We have analyzed the secular effects in the motion of an asteroid with negligible mass in inner case of a spatial restricted elliptic three-body problem with arbitrary inclination. Linear stability of the asteroid's orbits is studied. It is shown that equilibria, periodic orbits and other orbits of the double averaged planar restricted elliptic three-body problem are stable in the linear approximation as those orbits of the double averaged spatial restricted elliptic three-body problem with arbitrary inclination for all values of parameters. Numerical simulations of different orbits coincide with our results as well. This model with large eccentricities of the planet is of particular interest in relation to the study of the motion of exoplanets.

\section*{Acknowledgements}

The authors express their gratitude to Prof. Anatoly Neishtadt for suggestions on some of the topics in this work and to Prof. Xijun Hu and Prof. Yuwei Ou for discussions. Kaicheng Sheng thanks the National Natural Science Foundation of China (NSFC) for the support of this research (Grant:  12371192 \& 12271300). The work of Yiru Ye was supported by the XJTLU Research Development Funding (Grant: RDF-21-01-031).






\begin{thebibliography} {99}


\bibitem[\protect\citeauthoryear{Aksenov}{1979}]{akse}  Aksenov, E.P.: The doubly averaged, elliptical, restricted, three-body problem. Sov. Astronomy {\bf 23}, 236-240 (1979)
{
\bibitem[\protect\citeauthoryear{Arnold}{1961}]{Arnold1961}
Arnold, V.I.: The stability of the equilibrium position of a Hamiltonian system
of ordinary differential equations in the general elliptic case. Sov. Math., Dokl. {\bf 2}, 247–249 (1961)
}

\bibitem[\protect\citeauthoryear{Arnold et al.}{2006}]{akn} 
Arnold, V.I.,~Kozlov, V.V., Neishtadt, A.I.:
Mathematical Aspects of Classical and Celestial Mechanics, 3rd edn. Springer, New York (2006)

\bibitem[\protect\citeauthoryear{Brouwer \& Clemence}{1961}]{BC} Brouwer, D., Clemence, G.M.: Methods of Celestial Mechanics.  Academic Press, New York (1961)

\bibitem[\protect\citeauthoryear{Harrington}{1968}]{Harrington1968} Harrington, R.S.: {Dynamical evolution of triple stars.} Astron. J. {\bf 73}, 190-194 (1968) 


\bibitem[\protect\citeauthoryear{Katz \& Dong \& Malhotra}{2011}]{Katz2011} Katz, B., Dong, S., Malhotra, R.: {Long-term cycling of Kozai–Lidov cycles: extreme eccentricities and inclinations excited by a distant eccentric perturber.} Phys. Rev. Lett. {\bf 107}, 181101 (2011)

\bibitem[\protect\citeauthoryear{Kozai}{1962}]{kozai} Kozai, Y.: {Secular perturbations of asteroids with high inclination and eccentricity}.  Astron. J. {\bf 67}. 591-598 (1962)

\bibitem[\protect\citeauthoryear{Lei}{2022}]{Lei} Lei, H. . A Systematic Study about Orbit Flips of Test Particles Caused by Eccentric Von Zeipel–Lidov–Kozai Effects. The Astronomical Journal, 163 (2022).

\bibitem[\protect\citeauthoryear{Leontovich}{1962}]{Leontovich}Leontovich, A.M.: On the stability of the Lagrange periodic solutions of the restricted problem of three bodies. Sov. Math. Dokl. {\bf 3}, 425–428 (1962)

\bibitem[\protect\citeauthoryear{Lidov}{1962}]{lidov} Lidov, M.L.: {The evolution of orbits of artificial satellites of planets under the action of gravitational perturbations of external bodies}. Planet. Space Sci. {\bf 9}. 719-759 (1962)

\bibitem[\protect\citeauthoryear{Lidov \& Ziglin}{1974}]{lidov-ziglin} Lidov, M.L., Ziglin, S.L.:  The analysis of restricted circular twice-averaged
three body problem in the case of close orbits. Celest. Mech. {\bf 9}, 151-173 (1974)

\bibitem[\protect\citeauthoryear{Lithwick \& Naoz}{2011}]{Lithwick2011} Lithwick, Y., Naoz, S.: {The eccentric Kozai mechanism for a test particle.} Astrophys. J. 742:94 (2011)

\bibitem[\protect\citeauthoryear{Markeev}{1968}]{Markeev} Markeev, A.P.: Stability of a canonical system with two degrees of freedom in the presence of resonance. Journal of Applied Mathematics and Mechanics. {\bf 41}. 225-235 (1977)

\bibitem[\protect\citeauthoryear{Michtchenko \& Malhotra}{2004}]{Mi2004} Michtchenko, T.F., Malhotra, R.: Secular dynamics of the three-body problem: application to the $\upsilon$ Andromedae planetary system. ICARUS. {\bf 168}. 237-248 (2004)

\bibitem[\protect\citeauthoryear{Moiseev}{1945}]{moiseev} Moiseev, N.D.: {On some fundamental simplified schemes of celestial mechanics obtained by averaging of the restricted three-points problem. 2. On the averaged versions of the  three-dimensional  restricted circular three-points problem.} Trudy GAISh {\bf 15}.  100-117 (1945) (in Russian) 

\bibitem[\protect\citeauthoryear{Moser}{1968}]{Moser1968}
Moser, J.: Lectures on Hamiltonian Systems. Mem. Am. Math. Soc. {\bf 81}. American Mathematical Society, Providence, R. I. (1968)

\bibitem[\protect\citeauthoryear{Naoz}{2016}]{Naoz2016} Naoz, S.: {The eccentric Kozai–Lidov effect and its applications.} Ann. Rev. Astron. Astrophys. 54, 441–489 (2016)

\bibitem[\protect\citeauthoryear{Naoz et al.}{2013}]{naoz2013} Naoz, S., Farr, W. M., Lithwick, Y., Rasio, F. A., Teyssandier, J. (2013). Secular dynamics in hierarchical three-body systems. Monthly Notices of the Royal Astronomical Society, 431(3), 2155-2171.

\bibitem[\protect\citeauthoryear{Neishtadt}{1975}]{neish} Neishtadt, A.I.: {Stability of plane solutions in the doubly averaged restricted circular three-body problem}. Soviet Astronomy Letters {\bf 1}. 211-213 (1975)

\bibitem[\protect\citeauthoryear{Neishtadt \& Sheng  \& Sidorenko}{2021}]{NSS} Neishtadt, A.I., Sheng, K., Sidorenko, V.V.: {Stability analysis of apsidal alignment in double-averaged restricted elliptic three-body problem}. Celest. Mech. Dyn. Astr. 133:45 (2021)

\bibitem[\protect\citeauthoryear{Palaci\'an et al.}{2006}]{Palacian2006}
Palaci\'an, J.F., Yanguas, P., Fern\'andez, S., Nicotra, M.A.: {Searching for periodic orbits of the spatial elliptic restricted three-body problem by double averaging.} Physica D. {\bf 213}. 15-24 (2006) 


\bibitem[\protect\citeauthoryear{Shevchenko}{2016}]{shev} Shevchenko, I.: The Lidov-Kozai Effect - Applications in Exoplanet Research and Dynamical Astronomy. Springer, Berlin (2016)

\bibitem[\protect\citeauthoryear{Sidorenko}{2018}]{Sidorenko2018} Sidorenko, V.V.: {The eccentric Kozai-Lidov effect as a resonance phenomenon.} Celest. Mech. Dyn. Astr. 134:4 (2018)

\bibitem[\protect\citeauthoryear{Sokolskii}{1974}]{Sokolskii} Sokolskii, A.G.: {On the stability of an autonomous Hamiltonian system with two degrees of freedom in the case of equal frequencies.} Journal of Applied Mathematics and Mechanics, {\bf 38}. 741-749 (1974)

\bibitem[\protect\citeauthoryear{Stephen et al.}{2012}]{skdc} Stephen, R.K., David, R.C., Dawn, M.G., Kaspar, V.B.: {The exoplanet eccentricity distribution from Kepler planet candidates}. MNRAS {\bf 425}. 757 - 762 (2012)

\bibitem[\protect\citeauthoryear{Subaru Telescope Team}{2017}]{subaru}  Subaru Telescope Team: {Inclined Orbits Prevail in Exoplanetary Systems}.\\$https://subarutelescope.org/old/Pressrelease/2010/12/20/index.html$ (2017)


\bibitem[\protect\citeauthoryear{Vashkovyak}{1981}]{vashkovyak} Vashkovyak, M.A.:
Evolution of orbits in the restricted circular twice-averaged-three body problem. I. Qualitative
investigations. Cosm. Res. {\bf 19}. 1–10 (1981)

\bibitem[\protect\citeauthoryear{Vashkovyak}{1982}]{vash_planar}  Vashkovyak, M A.:  {Evolution of orbits in the two-dimensional restricted elliptical twice-averaged three-body problem}. Cosmic Research, {\bf 20}. 236-244 (1982)


\bibitem[\protect\citeauthoryear{von Zeipel}{1910}]{zeipel} von Zeipel, H.:  {Sur l'application des s\'eries de M. Lindstedt \`a l'\'etude du mouvement des com\'etes p\`eriodiques}. Astronomische Nachrichten {\bf 1983}. 345 - 418 (1910)

\bibitem[\protect\citeauthoryear{Ziglin}{1975}]{zig} Ziglin, S.L.: {Secular evolution of the orbit of a planet in a binary-star system}. Sov. Astron. Lett. {\bf 1}. 194-195 (1975)

\bibitem[\protect\citeauthoryear{Ziglin \& Lidov}{1977}]{Zig} Ziglin, S.L, Lidov, M.L.: Hill's case of the averaged problem of three bodies and the stability of plane orbits. Journal of Applied Mathematics and Mechanics {\bf 41}. 225-235 (1977)





\end{thebibliography}
\end{document}